\documentclass[11pt]{article}
\usepackage{latexsym}
\usepackage{amssymb}
\newcommand{\be}{\begin{equation}}
\newcommand{\ee}{\end{equation}\noindent}
\newcommand{\bear}{\begin{eqnarray}}
\newcommand{\ear}{\end{eqnarray}\noindent}

\date{}
\renewcommand{\theequation}{\arabic{equation}}

\newcommand{\slD}{\raise.15ex\hbox{$/$}\kern-.57em\hbox{$D$}}
\newcommand{\slpartial}{\raise.15ex\hbox{$/$}\kern-.57em\hbox{$\partial$}}
\newcommand{\slG}{{{\dot G}\!\!\!\! \raise.15ex\hbox {/}}}

\usepackage[colorlinks,urlcolor=blue]{hyperref}

\def\hhref#1{\href{http://arxiv.org/abs/#1}{arXiv:#1}} 

\begin{document}
\pagestyle{empty}
\renewcommand{\thefootnote}{\fnsymbol{footnote}}
\hfill {\sl}

\hfill {\sl}   
\vskip-.1pt
\hfill {\large IHES/P/04/31 (note added, January 2013)}
\vskip .4cm
\begin{center}
{\Large\bf Bernoulli Number Identities from Quantum Field Theory and Topological String Theory}
\vskip1.3cm
{\large Gerald V. Dunne}
\\[1.5ex]
{\it
Department of Physics and Department of Mathematics\\
University of Connecticut\\
Storrs CT 06269, USA
}
\vspace{.8cm}

 {\large Christian Schubert
}
\\[1.5ex]
{\it
Institut des Hautes \'Etudes Scientifiques
\\
Le Bois-Marie, 35, F-91440 Bures-sur-Yvette, FRANCE
\\
\vspace{10pt}
{\rm and}
\vspace{10pt}
\\
Department of Physics and Geology
\\
University of Texas Pan American
\\
Edinburg, TX 78541-2999, USA
\\
\vspace{10pt}
{\rm and}
\vspace{10pt}
\\
Instituto de F\'{\i}sica y Matem\'aticas
\\
Universidad Michoacana de San Nicol\'as de Hidalgo\\
Apdo. Postal 2-82\\
C.P. 58040, Morelia, Michoac\'an, M\'exico\\
}
\vspace{1cm}

\vspace{.2cm}
 {\large \bf Abstract}
\end{center}
\begin{quotation}
\noindent
We present a new method for the derivation of
convolution identities for finite sums of products of Bernoulli numbers. 
Our approach is motivated by the role of 
these identities in quantum field theory and string theory.
We first show that the Miki identity and 
the Faber-Pandharipande-Zagier (FPZ) identity are closely related, 
and give 
simple unified proofs which naturally yield a new Bernoulli number 
convolution 
identity. We then generalize each of these three identities into new 
families of convolution identities depending on a continuous parameter.
We rederive a cubic generalization of Miki's identity due to Gessel
and obtain a new similar identity generalizing the FPZ identity.
The generalization of the method to the derivation of
convolution identities of arbitrary order is outlined.
We also describe an extension to identities which relate convolutions of Euler and Bernoulli
numbers.

\end{quotation}
\vskip 1cm
\clearpage
\renewcommand{\thefootnote}{\protect\arabic{footnote}}
\pagestyle{plain}

\section{Introduction: Convolution identities for Bernoulli numbers}
\label{introduction}
\renewcommand{\theequation}{1.\arabic{equation}}
\setcounter{equation}{0}

The Bernoulli numbers $B_n$ are defined by the generating function \cite{abramowitz}
\bear
b(x) \equiv {x\over e^x-1} = \sum_{n=0}^{\infty}B_n{x^n\over n!}
\label{defbernoulli}
\ear
As is well known, the  $B_n$ play an important role in combinatorics 
and 
number theory, 
and there exist many combinatorial identities involving these numbers 
\cite{comtet,apostol}.
In the present paper, we will be concerned with the special case
of convolution identities, which involve finite sums of products of 
Bernoulli 
numbers. The best known such convolution identity was found already by 
Euler 
(and independently by Ramanujan):

\underline{{\sl Theorem 1.1 :}} (Euler; Ramanujan) For integer $n\geq 2$,
\bear
\sum_{k=1}^{n-1}{2n\choose 2k}B_{2k}B_{2n-2k}=-(2n+1)B_{2n}
\label{t11}
\ear
The proof follows directly from the definition (\ref{defbernoulli}), 
noting 
that the generating
function $b(x)$ satisfies
\bear
b(x)^2 = (1-x)b(x)-xb'(x)
\label{idb}
\ear
Equivalently, (\ref{t11}) follows by comparing the series expansions 
of either 
side of the trigonometric identity
\bear
\coth^2 x=1-(\coth x)^\prime
\ear

Many more such identities, involving folded sums of Bernoulli numbers, 
have 
been
found since Euler's work (see, e.g., \cite{dilcher,huahua,basu}).  
Most of 
them are similar
to Euler's identity (\ref{t11}) in the sense that they involve 
$B_n\over n!$
rather than $B_n$ itself, as is the case already for the defining 
formula (\ref
{defbernoulli}).
Identities involving the $B_n$ themselves, without the factorial 
denominator,
are much rarer. One such identity was found by H. Miki in 1978.

\underline{{\sl Theorem 1.2 :}} (Miki \cite{miki}): For integer $n\geq 
2$,
\bear
\sum_{k=1}^{n-1}{B_{2k}B_{2n-2k}\over
  (2k)(2n-2k)}=\sum_{k=1}^{n-1}{B_{2k}B_{2n-2k}\over
  (2k)(2n-2k)}{2n\choose 2k} +{B_{2n}\over n}
H_{2n}
\label{mikiidentity}
\ear
Here $H_i$ denotes the $i$th harmonic number,

\bear
H_i &\equiv& \sum_{j=1}^i {1\over j}
\label{defH}
\ear
As is well-known, the harmonic numbers can be alternatively expressed in terms of the
digamma function $\psi(x)=\Gamma'(x)/\Gamma(x)$ 
and the Euler- Mascheroni constant $\gamma$:

\bear
H_i &=& \psi(i+1)+\gamma
\label{Htopsi}
\ear
Miki's original proof of  the identity (\ref{mikiidentity}) 
identity is quite involved \cite{miki}. A more
elementary proof was given recently by I. M. Gessel \cite{gessel}, 
using two 
different expressions for the
Stirling numbers of the second kind. In Section 2.1 we present an even 
simpler 
proof, based on an appropriate generating function.

In 1998 C. Faber and R. Pandharipande \cite{fabpan} found that
certain conjectural relations between Hodge integrals in Gromov-Witten 
theory 
\cite{eghoxi,getpan} (see also \cite{lilizh})
require the following identity to hold.

\underline{{\sl Theorem 1.3 :}} (Faber and Pandharipande, with a proof 
by 
Zagier \cite{fabpan}) 
For integer $n\geq 2$,
\bear
\sum_{k=1}^{n-1}{\bar{B}_{2k}\bar{B}_{2n-2k}\over
  (2k)(2n-2k)}
 &=&\frac{1}{n}\sum_{k=1}^{n}{B_{2k}\bar{B}_{2n-2k}\over
  (2k)}{2n\choose 2k}  +{\bar{B}_{2n}\over n}H_{2n-1}
\label{fpzidentity}\\
{\rm where}\quad   \bar{B}_n&\equiv& \Bigl({1-2^{n-1}\over 
2^{n-1}}\Bigr)B_n
\label{bbar}
\ear
A proof of the Faber-Pandharipande-Zagier (FPZ) identity 
(\ref{fpzidentity}) 
was given by D. Zagier in an appendix to \cite{fabpan}.
Note that the structure of the FPZ identity (\ref{fpzidentity}) is 
similar to 
Miki's identity (\ref{mikiidentity}). We later show (see Theorem 2.1) 
that this 
similarity is even more striking if Miki's identity is written in a 
slightly 
different form.

Apart from pure mathematics, the Bernoulli numbers
appear prominently in perturbative quantum field theory. This comes 
about 
at a very basic level \cite{ss1}: perturbative loop calculations in
quantum field theory generally involve traces of inverse powers of 
derivatives
of functions defined on a circle. Since
the spectrum of the ordinary derivative operator $\partial_P$ with 
periodic 
boundary conditions
consists of the integer numbers, one has
\bear
{\rm tr}(\partial_P^{-2n}) \sim \sum_{k=1}^{\infty}{1\over k^{2n}} = 
\zeta(2n)
\label{invtrace}
\ear
But $\zeta(2n)$ is related to the Bernoulli numbers through Euler's 
identity,
  \bear
 B_{2n} &=& (-1)^{n+1}2{(2n)!\over (2\pi)^{2n}}\zeta(2n)
 \label{ideulerzeta}
\ear
The Bernoulli numbers also appear naturally in so-called "effective 
action" 
computations in quantum field theory, a field pioneered by W. 
Heisenberg and H. 
Euler, and V. Weisskopf \cite{dunne}.
 
In \cite{sd1,sd2}, the present authors found that Miki's identity 
arises 
naturally in a certain computation in perturbative quantum field 
theory. 
Specifically, it arises in the course of the  calculation of the 
two-loop 
effective
Lagrangian for quantum electrodynamics in a constant background 
self-dual 
field. This calculation was       
done using two different integral representations. It turned out that 
both
representations yield a result for the coefficients of the weak field 
expansion of this
effective Lagrangian which involve a convolution of Bernoulli numbers, 
and 
that 
what is needed to show the equivalence of both results is precisely 
Miki's 
identity 
(\ref{mikiidentity}). Since the two integral representations used are 
related 
by a
simple coordinate transformation,
this actually yields a new, and quite straightforward, proof of Miki's 
identity. 
This proof will be given in section \ref{mikiproof}. 
The simplicity of the approach presented here suggests a number of 
generalizations of these
two identities, some of which are presented in sections 
\ref{newidentity}
to \ref{multi}. It is straightforward to verify these identities explicitly using symbolic computer programs, and we have performed such checks for these identities for various values of $n$.
 Further  generalizations  are
outlined in the conclusions.

\section{Simple generating function proofs of Miki and FPZ identities}
\label{mikifpz}
\renewcommand{\theequation}{2.\arabic{equation}}
\setcounter{equation}{0}
In this section we present simple unified proofs of the Miki and FPZ 
identities based on generating functions.

\subsection{Generating function proof of Miki's identity}
\label{mikiproof}
Consider the generating function
\bear
\tilde \psi(x) &\equiv & \psi(x) - \ln x + {1\over 2x}
\label{defpsitilde}
\ear
This function plays an important role in the quantum field theory computations in \cite{sd1,sd2}.
From the asymptotic (large x) expansion of the digamma function \cite
{abramowitz}
it follows that
\bear
\tilde \psi(x)\, \sim  - \sum_{k=1}^{\infty} {B_{2k}\over 
2k}\frac{1}{x^{2k}}
\label{asymptildepsi}
\ear
Thus for the square of $\tilde\psi$ one finds
\bear
\bigl[\tilde\psi (x)\bigr]^2 \sim 
\sum_{n=2}^{\infty} {1\over x^{2n}}
\sum_{k=1}^{n-1}{B_{2k}B_{2n-2k}\over (2k)(2n-2k)}
\label{asymptildepsisquared}
\ear
Thus, we see that the $[\tilde{\psi}(x)]^2$ is the generating function 
for the 
left-hand side of Miki's identity (\ref{mikiidentity}). We prove Miki's 
identity by comparing (\ref{asymptildepsisquared}) with the square of 
the 
following integral representation [see Eq 1.7.2 (25) in \cite{bateman}] of $\tilde{\psi}(x)$ :
\bear
\tilde \psi(x) &=& -\int_0^{\infty}ds\,{\rm e}^{-2xs}
\,\Bigl(\coth s - {1\over s}\Bigr)
\label{intreppsitilde}
\ear
 We break this comparison into three straightforward lemmas.

\underline{{\sl Lemma 2.1.1 :}} 
\bear 
\left[\tilde{\psi}(x)\right]^2&=& \int_0^\infty dy\, y\, e^{-2 x 
y}\int_0^1 du 
\,\left\{-1+2\left(\coth y-\frac{1}{y}\right)\left(\coth 
yu-\frac{1}{yu}\right) 
\right.
\nonumber\\
&& \left. -\frac{2}{y(1-u)}\left[u\left( \coth yu-\frac{1}{yu}\right)- 
\left( 
\coth y-\frac{1}{y}\right)\right] \right\}
\label{l211}
\ear
\underline{{\sl Proof :}} Squaring (\ref{intreppsitilde}) we find
\bear
\bigl[\tilde\psi(x)\bigr]^2 =\int_0^{\infty}ds\int_0^{\infty}ds' 
\,{\rm e}^{-2x
(s+s')}
\Bigl\{\coth s\coth s' -\Bigl({1\over s}\coth s' + {1\over s'}\coth s 
\Bigr) 
+{1\over ss'}\Bigr\} 
\label{trafo2}
\ear
Now, using the trigonometric identity
\bear
\coth s \coth s' &=& \coth (s+s') \Bigl(\coth s + \coth s' \Bigr) -1 
\quad ,
\label{trigcoths}
\ear
together with the symmetry $s\leftrightarrow s'$, 
and the transformation of variables
\footnote{In quantum field theory terms, this change of variables 
corresponds to a change from a Feynman parameter integral (see, e.g., 
\cite
{itzzub})
to a worldline parameter integral \cite{ss2}.} 
\bear
y = s+s^\prime \qquad , \qquad u = {s'\over s+s'}\quad ,
\label{trafofeynwl}
\ear
it is straightforward to show that $[\tilde{\psi}(x)]^2$  can be 
rewritten as 
in {\sl Lemma 2.1.1}. Note that the change of variables 
(\ref{trafofeynwl}) 
introduces a Jacobian factor of $y$.   $\square$

Our proof of Miki's identity follows by evaluating the asymptotic 
expansions 
of the integrals appearing on the right-hand side of (\ref{l211}).

\underline{{\sl Lemma 2.1.2 :}} 
\bear
2\int_0^{\infty}dy \,y \,{\rm e}^{-2xy} \int_0^1 du
\Bigl(\coth y - {1\over y}\Bigr)
\Bigl(\coth yu - {1\over yu}\Bigr)
\sim
\sum_{n=2}^{\infty}{1\over x^{2n}}
\sum_{k=1}^{n-1} 
{B_{2k}B_{2n-2k}\over (2k)(2n-2k)}{2n\choose 2k}\nonumber\\
\label{l212}
\ear
\underline{\sl Proof :} The $u$ integral is elementary:
\bear
\int_0^1 du \Bigl( \coth yu - {1\over yu}\Bigr) &=&  {1\over y}{\rm 
ln}\, \Bigl
({\sinh y\over y}\Bigr)
\label{intcothyu}
\ear
After an integration by parts, the remaining $y$ integral takes the 
form
\bear
2x\int_0^{\infty}dy \,{\rm e}^{-2xy} \,
{\rm ln}^2\, \Bigl({\sinh y\over y}\Bigr)\nonumber\\
\label{firstint}
\ear
The asymptotic expansion of the $y$ integral is obtained using the 
Taylor 
expansion \cite{abramowitz}
\bear
\ln \Bigl({\sinh y\over y}\Bigr) = 
\sum_{k=1}^{\infty}{2^{2k-1}B_{2k}\over k
(2k)!} y^{2k}
\label{taylsinh}
\ear
which directly yields the result (\ref{l212}) after performing the $y$ 
integration. $\square$

\underline{{\sl Lemma 2.1.3 :}} 
\bear
&&-2\int_0^\infty dy\, y\, e^{-2 x y}\int_0^1 du \, 
\left\{\frac{1}{y(1-u)}
\left[u\left( \coth yu-\frac{1}{yu}\right)- \left( \coth 
y-\frac{1}{y}\right)
\right] \right\}\nonumber\\
&& \qquad \qquad \sim \quad  \sum_{n=1}^\infty \frac{1}{x^{2n}} 
\frac{B_{2n}}
{n} H_{2n}
\label{l213}
\ear
\underline{\sl Proof :} First, consider the $u$ integral. We use the 
Taylor 
expansion of the $\coth$ function \cite{abramowitz},
\bear
y\coth y = \sum_{k=0}^{\infty}{2^{2k}B_{2k}
\over (2k)!} y^{2k}
\label{taylcoth}
\ear
for both $\coth yu$ and $\coth y$. The $u$ integral becomes 
elementary:
\bear
\int_0^1 du \frac{\left[u\left( \coth yu-\frac{1}{yu}\right)- \left( 
\coth y-
\frac{1}{y}\right)\right]}{(1-u)}
& = & \sum_{n=1}^\infty \frac{B_{2n} 2^{2n} y^{2n-1}}{(2n)!} \int_0^1 
du   
\left(\frac{u^{2n}-1}{1-u}\right)\nonumber\\
&= &-\sum_{n=1}^\infty \frac{B_{2n} 2^{2n} y^{2n-1}}{(2n)!} 
H_{2n}\nonumber\\
\ear
Doing the $y$ integral we obtain the result of {\sl Lemma 2.1.3}. $\square$

Miki's identity (\ref{mikiidentity}) is then proved by comparing the 
results 
of {\sl Lemmas} 2.1.1 - 2.1.3 with (\ref{asymptildepsisquared}).

We conclude this section on Miki's identity by remarking that in the 
proof of 
{\sl Lemma 2.1.2}, the partial integration in $y$ leading to 
(\ref{firstint}) 
is not essential. If, instead, one does the $y$  integral directly 
using (\ref
{taylcoth})
and (\ref{taylsinh}), one arrives at a slightly different version of 
Miki's 
identity:

\underline{\sl Theorem 2.1 :} (Modified form of Miki's identity): For 
integer 
$n\geq 2$, 
\bear
\sum_{k=1}^{n-1}{B_{2k}B_{2n-2k}\over
  (2k)(2n-2k)}&=&{1\over n}\sum_{k=1}^{n-1}{B_{2k}B_{2n-2k}\over
  (2k)}{2n\choose 2k} +{B_{2n}\over n}H_{2n} \nonumber\\
&=&{1\over n}\sum_{k=1}^n{B_{2k}B_{2n-2k}\over
  (2k)}{2n\choose 2k} +{B_{2n}\over n}H_{2n-1}
\label{t21}
\ear
where we have used $H_{2n}=H_{2n-1}+{1\over 2n}$, and $B_0=1$.

\underline{\sl Comment 2.1.1 :} This last form (\ref{t21}) of Miki's 
identity 
brings out most clearly the similarity to the FPZ identity 
(\ref{fpzidentity}).

\subsection{Generating function proof of the FPZ identity}
\label{fpzproof}
To prove the FPZ identity we use, instead of $\tilde\psi(x)$, the 
generating 
function $\bar{\psi}(x)$ defined by
\bear
\bar \psi(x) &\equiv & \psi(x+\frac{1}{2}) - \ln x 
\label{defpsibar}
\ear
The large $x$ expansion of $\bar{\psi}(x)$ is
\bear
\bar\psi(x) &\sim& - \sum_{k=1}^{\infty}{\bar B_{2k}\over 
2k}\frac{1}{x^{2k}}
\label{psibarlargex}
\ear
where $\bar{B}_{2k}$ was defined in (\ref{bbar}). 
The expansion (\ref{psibarlargex}) follows from the corresponding 
expansion 
(\ref{asymptildepsi}) for $\tilde{\psi}(x)$, 
using the ``doubling'' identity \cite{abramowitz} for the $\psi$ 
function,
\bear
\psi(2x) &=&  {1\over 2} \psi(x) + {1\over 2} \psi(x+{1\over 2}) + \ln 
2
\label{iddouble}
\ear
Thus, the square of $\bar\psi(x)$ is the generating function for the 
left-hand 
side of the FPZ identity (\ref{fpzidentity}) :
\bear
\bigl[\bar\psi (x)\bigr]^2 \sim 
\sum_{n=2}^{\infty} {1\over x^{2n}}
\sum_{k=1}^{n-1}{\bar B_{2k}\bar B_{2n-2k}\over (2k)(2n-2k)}
\label{asympbarpsisquared}
\ear
The generating function for the right-hand side of the FPZ identity is 
obtained by squaring the following integral representation for 
$\bar\psi(x)$:
\bear
\bar\psi(x) = - \int_0^{\infty} ds\, {\rm e}^{-2xs}
\Bigl ({1\over\sinh s} -{1\over s}
\Bigr)
\label{intreppsibar}
\ear

\underline{\sl Lemma 2.2.1 :}
\bear 
\left[\bar{\psi}(x)\right]^2&=& 2\int_0^\infty dy\, y\, e^{-2 x 
y}\int_0^1 du 
\,\left\{\frac{1}{\sinh y}\left(\coth yu-\frac{1}{yu}\right) \right.
\nonumber\\
&& \left. -\frac{1}{y(1-u)}\left[u\left( \frac{1}{\sinh 
yu}-\frac{1}{yu}
\right)- \left( \frac{1}{\sinh y}-\frac{1}{y}\right)\right] \right\}
\label{l221}
\ear
\underline{\sl Proof :} Squaring the integral representation (\ref
{intreppsibar}), using the trigonometric identity
\bear
{1\over \sinh(s)\sinh(s')} &=& 
{\coth(s)+\coth(s')\over\sinh(s+s')}\quad ,
\label{trigsinh}
\ear
and the symmetry under $s\leftrightarrow s^\prime$, it follows that
\bear
\bigl[\bar\psi(x)\bigr]^2 &=&
2\int_0^{\infty}ds\int_0^{\infty}\,ds' \,{\rm e}^{-2(s+s')x}
\biggl\lbrace
{1\over\sinh(s+s')}\Bigl(\coth(s')-{1\over s'}\Bigr)
\nonumber\\&&\hspace{40pt}
+{1\over s}\Bigl({1\over\sinh(s+s')}-{1\over 
s+s'}-{1\over\sinh(s')}+{1\over 
s'}\Bigr)
\biggr\rbrace
\label{trafo2bar}
\ear
Applying the transformation of variables (\ref{trafofeynwl}), we 
obtain (\ref{l221}). $\square$

Our proof of the FPZ identity now follows by evaluating the asymptotic 
expansions of the integrals appearing on the right-hand side of 
(\ref{l221}).

\underline{\sl Lemma 2.2.2 :}
\bear 
2\int_0^\infty dy\, y\, e^{-2 x y}\int_0^1 du \, \left\{\frac{1}{\sinh 
y}\left
(\coth yu-\frac{1}{yu}\right) \right\} \sim   \sum_{n=1}^\infty 
\frac{1}{x^
{2n}} \frac{1}{n} \sum_{k=1}^{n} \frac{B_{2k} 
\bar{B}_{2n-2k}}{(2k)}{2n\choose 
2k}
\label{l222}
\ear
\underline{\sl Proof :} The proof is almost identical to the proof of 
{\sl 
Lemma 2.1.2}, but in doing the $y$ integral we use an asymptotic 
expansion of 
$1/\sinh y$ rather than $\coth y$. This has the effect of replacing one 
of the Bernoulli number factors $B_{2n-2k}$ by $\bar{B}_{2n-2k}$, and also of changing the upper 
limit of 
the $k$ summation from $(n-1)$ to $n$. $\square$

\underline{\sl Lemma 2.2.3 :}
\bear 
&&- 2\int_0^\infty dy\, y\, e^{-2 x y}\int_0^1 du \, 
\left\{\frac{1}{y(1-u)}
\left[u\left( \frac{1}{\sinh yu}-\frac{1}{yu}\right)- \left( 
\frac{1}{\sinh y}-
\frac{1}{y}\right)\right] \right\}\nonumber\\
&&\quad \sim \quad  \sum_{n=1}^\infty \frac{1}{x^{2n}} 
\frac{\bar{B}_{2n}}{n} 
H_{2n-1}
\label{l223}
\ear
\underline{\sl Proof :} The proof is almost identical to the proof of 
{\sl 
Lemma 2.1.3}, except we use the asymptotic expansion of $1/\sinh y$ 
rather than 
$(\coth y-\frac{1}{y})$, which has the effect of replacing $B_{2n}$ by 
$\bar{B}_{2n}$.  $\square$

The FPZ identity (\ref{fpzidentity}) is then proved by comparing the 
results 
of {\sl Lemmas} 2.2.1 - 2.2.3 with (\ref{asympbarpsisquared}). Note 
that the 
$n=1$ terms on the right-hand sides of (\ref{l222}) and (\ref{l223}) 
cancel, permitting the comparison with (\ref{asympbarpsisquared}).

 \section{A new convolution identity}
\label{newidentity}
\renewcommand{\theequation}{3.\arabic{equation}}
\setcounter{equation}{0}

The similarity between the proofs and forms of the Miki and FPZ 
identities 
immediately suggests a new identity, in which on the left-hand side the 
$B_{2n}
$ and $\bar{B}_{2n}$ are mixed. As is clear from Section 2, such an 
identity 
could be derived using the generating function $\tilde{\psi}(x)\,\bar{\psi}(x)$, and 
comparing its summation and integral representations. However, there is 
another, even simpler, way to derive this mixed identity.  Note  that 
the two 
generating functions $\tilde{\psi}(x)$ and $\bar{\psi}(x)$ are related via the $\psi$ function doubling 
identity (\ref{iddouble}) as:
\bear
 \tilde{\psi}(x)+ \bar{\psi}(x)=2 \tilde{\psi}(2x)
 \label{crossrelation}
\ear
Thus, it follows that
\bear
2\, \tilde{\psi}(x)\, \bar{\psi}(x)=4 \left[ \tilde{\psi}(2x)\right]^2 
- \left
[ \tilde{\psi}(x)\right]^2 - \left[ \bar{\psi}(x)\right]^2
\label{cross}
\ear
A new Bernoulli convolution identity emerges by using the asymptotic 
expansions (\ref{asymptildepsi}) and (\ref{psibarlargex}) for the 
left-hand 
side of (\ref{cross}), and the asymptotic expansions of the squares of 
the 
integral representations (\ref{intreppsitilde}) and (\ref{intreppsibar}) of $\tilde{\psi}(x)$ and 
$\bar{\psi}(x)
$, respectively, on the right-hand side. All necessary results for the 
squares 
of the relevant integral representations are contained in Lemmas 2.1.1 
- 2.1.3, 
and 2.2.1 - 2.2.
3. To express the result in a symmetrical form, we use the modified 
form of 
Miki's identity in (\ref{t21}).

\underline{\sl Theorem 3.1 :} For integer $n\geq 2$:
\bear
\sum_{k=1}^{n-1} \frac{B_{2k} \bar{B}_{2n-2k}}{(2k)(2n-2k)} = 
\frac{1}{n}\sum_{k=1}^{n}\frac{B_{2k} B_{2n-2k}}{(2k)}{2n\choose 2k} 
\Bigl
({1-2^{2k-1}\over 2^{2n-1}}\Bigr)+ \frac{1}{n}{B_{2n}\over 
2^{2n}}H_{2n-1}
\label{t31}
\ear
\underline{\sl Proof :} The generating function for the left-hand side 
is 
given by half the left-hand side of (\ref{cross}). The right-hand side 
is 
obtained by using Miki's identity in the form (\ref{t21}) for the 
squares of 
$\tilde{\psi}$, and the FPZ 
identity (\ref{fpzidentity}) for the square of $\bar{\psi}$. Simple 
algebra 
then leads to the form in (\ref{t31}).

\section{Three infinite families of convolution identities}
\label{families}
\renewcommand{\theequation}{4.\arabic{equation}}
\setcounter{equation}{0}

The use of the generating functions $\tilde{\psi}(x)$ and 
$\bar{\psi}(x)$ to 
prove the Miki and FPZ identities, as well as the new ``crossed" 
identity (\ref
{t31}) in {\sl Theorem 3.1}, immediately leads to natural 
generalizations of each type of identity.

\subsection{Generalization of Miki's Identity}
\label{mikigen}

To derive a generalization of Miki's identity, consider the $p^{\rm th}$ 
derivative of the generating function $\tilde{\psi}(x)$. This has the 
large $x$ 
asymptotic expansion:
\bear
\tilde{\psi}^{(p)}(x)\sim (-1)^{p+1}\sum_{n=1}^\infty 
\frac{B_{2n}\Gamma(2n+p)}
{(2n) \Gamma(2n)}\, \frac{1}{x^{2n+p}}
\label{asymptoticptder}
\ear
This function also has the following integral representation:
\bear
\tilde{\psi}^{(p)}(x)= -(-2)^{p}\int_0^{\infty}ds\,{\rm e}^{-2xs}
\,s^p\, \Bigl(\coth s - {1\over s}\Bigr)
\label{integralptder}
\ear
We can use this integral representation to extend the definition of $\tilde{\psi}^{(p)}(x)$ to non-integer
values of $p$; for this extrapolation the expansion (\ref{asymptoticptder}) continues to hold,
as can be seen by using (\ref{taylcoth}) under the integral in (\ref{integralptder}).
Thus, in the following let $p$ denote an arbitrary non-negative number. 
We can derive new identities, for any such $p$, by squaring 
these two representations of $\tilde{\psi}^{(p)}(x)$, and then 
comparing, just 
as was done (for $p=0$) to prove Miki's identity.
The proof proceeds in a very similar manner.

\underline{\sl Lemma 4.1.1 :}
\bear
\left[\tilde{\psi}^{(p)}(x)\right]^2 \sim \sum_{n=2}^\infty 
\frac{1}{x^
{2n+2p}} \sum_{k=1}^{n-1} \frac{B_{2k} 
B_{2n-2k}}{(2k)(2n-2k)}\frac{\Gamma(2k+p)
\Gamma(2n-2k+p)}{ \Gamma(2k)\Gamma(2n-2k)}\, 
\label{l411}
\ear
\underline{\sl Proof:} follows from (\ref{asymptoticptder}). $\square$

\underline{\sl Lemma 4.1.2 :}
\bear 
\left[\tilde{\psi}^{(p)}(x)\right]^2&=& 2^{2p}\int_0^\infty dy\, 
y^{2p+1}\, e^
{-2 x y}\int_0^1 du \, u^p(1-u)^p\left\{-1+2\left(\coth 
y-\frac{1}{y}\right)
\left(\coth yu-\frac{1}{yu}\right) \right.
\nonumber\\
&& \left. -\frac{2}{y(1-u)}\left[u\left( \coth yu-\frac{1}{yu}\right)- 
\left( 
\coth y-\frac{1}{y}\right)\right] \right\}
\label{l412}
\ear
\underline{\sl Proof :} Square the integral representation (\ref
{integralptder}), change variables from $s$ and $s^\prime$ to $y$ and 
$u$, as 
in (\ref{trafofeynwl}), and regroup terms as in the proof of {\sl Lemma 
2.1.1}.
Note that the argument about
symmetrizing with respect to $s$ and $s^\prime$ still holds because 
the extra 
factors in the integrand appear as $(s \, 
s^\prime)^p=y^{2p}u^p(1-u)^p$. $\square$

Now consider each of the three terms appearing on the RHS of 
(\ref{l412}).

\underline{\sl Lemma 4.1.3 :} 
\bear
-2^{2p}\int_0^\infty dy\, y^{2p+1}\, e^{-2 x y}\int_0^1 du \, 
u^p(1-u)^p \sim -
\frac{\Gamma^2(p+1)}{4 x^{2p+2}}
\label{l413}
\ear
\underline{\sl Proof :} immediate. $\square$

\underline{\sl Lemma 4.1.4 :} 
\bear 
&&2^{2p+1}\int_0^\infty dy\, y^{2p+1}\, e^{-2 x y}\int_0^1 du \, 
u^p(1-u)
^p\left(\coth y-\frac{1}{y}\right)\left(\coth yu-\frac{1}{yu}\right) 
\nonumber\\
&& \hskip 1cm \sim \quad 2 \Gamma(p+1)\sum_{n=2}^\infty \frac{1}{x^{2n+2p}} 
\sum_{k=1}^
{n-1} \frac{B_{2k} 
B_{2n-2k}}{(2k)!(2n-2k)!}\frac{\Gamma(2k+p)\Gamma(2n+2p)}
{\Gamma(2p+2k+1)}
\label{l414}
\ear
\underline{\sl Proof :} First, consider the $u$ integral:
\bear 
\int_0^1 du \, u^p(1-u)^p\left(\coth yu-\frac{1}{yu}\right)  &=& 
\sum_{n=1}
^\infty \frac{B_{2n} 2^{2n} y^{2n-1}}{(2n)!} \int_0^1 du \, 
u^{2n+p-1}(1-u)
^p\nonumber\\
&=& \sum_{n=1}^\infty \frac{B_{2n} 2^{2n} y^{2n-1}}{(2n)!} 
\frac{\Gamma(p+1)
\Gamma(p+2n)}{\Gamma(2p+2n+1)}
\ear
Now doing the $y$ integral we obtain:
\bear
&& 2^{2p+1}\sum_{n=1}^\infty \frac{B_{2n} 2^{2n}}{(2n)!} 
\frac{\Gamma(p+1)
\Gamma(p+2n)}{\Gamma(2p+2n+1)}  \int_0^\infty dy\, y^{2p+2n}\, e^{-2 x 
y}\left
(\coth y-\frac{1}{y}\right)\nonumber\\
&&\quad \sim \quad 2^{2p+1}\sum_{n=1}^\infty \frac{B_{2n} 
2^{2n}}{(2n)!} \frac
{\Gamma(p+1)\Gamma(p+2n)}{\Gamma(2p+2n+1)} \sum_{k=1}^\infty 
\frac{B_{2k} 2^
{2k}}{(2k)!} \frac{\Gamma(2p+2n+2k)}{(2x)^{2p+2n+2k}}
\ear
from which {\sl Lemma 4.1.4}  follows. $\square$

\underline{\sl Lemma 4.1.5 :}
\bear
&&-2^{2p+1}\int_0^\infty dy\, y^{2p+1}\, e^{-2 x y}\int_0^1 du \, 
u^p(1-u)
^p\left\{\frac{1}{y(1-u)}\left[u\left( \coth yu-\frac{1}{yu}\right)- 
\left( 
\coth y-\frac{1}{y}\right)\right] \right\}\nonumber\\
&& \quad \sim \quad 2 \sum_{n=2}^\infty \frac{1}{x^{2n+2p}} 
\frac{B_{2n} \Gamma
(2n+2p)}{(2n)!} \sum_{k=1}^{2n} \beta(p+k,p+1)
\label{l44}
\ear
\underline{\sl Proof :}
First, consider the $u$ integral:
\bear
&&\int_0^1 du \, u^p(1-u)^{p-1}\left[u\left( \coth 
yu-\frac{1}{yu}\right)- 
\left( \coth y-\frac{1}{y}\right)\right]\nonumber\\
&&\quad \sim \quad \sum_{n=1}^\infty \frac{B_{2n} 2^{2n} 
y^{2n-1}}{(2n)!} 
\int_0^1 du  \, u^p(1-u)^{p-1} \left(u^{2n}-1\right)\nonumber\\
&&\quad \sim \quad \sum_{n=1}^\infty \frac{B_{2n} 2^{2n} 
y^{2n-1}}{(2n)!} \sum_
{k=1}^{2n} \beta(p+k,p+1)
\ear
where $\beta(p,q)$ is the Euler beta function.
Doing the $y$ integral we obtain the result of {\sl Lemma} 4.1.5. $\square$

We are now ready to state the generalization of Miki's identity:

\underline{\sl Theorem 4.1 :} For any  $p\geq 0$, and for 
integer 
$n\geq 2$:
\bear
\sum_{k=1}^{n-1} \frac{B_{2k} 
B_{2n-2k}}{(2k)(2n-2k)}\frac{\Gamma(2k+p)\Gamma
(2n-2k+p)}{ \Gamma(2k)\Gamma(2n-2k)} &=& 2\Gamma(p+1)\sum_{k=1}^{n} 
\frac{B_{2k} B_{2n-
2k}}{(2k)!(2n-2k)!}\frac{\Gamma(2k+p)\Gamma(2n+2p)}{\Gamma(2p+2k+1)}\nonumber\\
&& +2  \frac{B_{2n} \Gamma(2n+2p)}{(2n)!} \sum_{k=1}^{2n-1} 
\beta(p+k,p+1)
\label{t41}
\ear

\underline{\sl Proof :} Follows by comparing the result of {\sl 
Lemma 4.1.1} 
with those of {\sl Lemmas} 4.1.2 - 4.1.5. $\square$ 

\underline{\sl Comment 4.1.1 :} When $p=0$ we recover from {\sl 
Theorem 4.1}  
Miki's identity in the form of {\sl Theorem 2.1}. 

\underline{\sl Comment 4.1.2 :} When $p=1$ we obtain from {\sl Theorem 
4.1} a 
convolution identity which just involves the Bernoulli numbers 
themselves on 
the left-hand side: for $n\geq 2$,
\bear
\sum_{k=1}^{n} B_{2k} B_{2n-2k}= \frac{1}{n+1} \sum_{k=1}^{n} B_{2k} 
B_{2n-2k}
{2n+2\choose 2k+2} + 2n B_{2n}
\label{c412}
\ear

\subsection{Generalization of the FPZ Identity}
\label{fpzgen}

To derive a generalization of the FPZ identity, consider the $p^{\rm th}$ 
derivative of the generating function $\bar{\psi}(x)$. This has the 
large $x$ 
asymptotic expansion:
\bear
\bar{\psi}^{(p)}(x)\sim (-1)^{p+1}\sum_{n=1}^\infty 
\frac{\bar{B}_{2n}\Gamma
(2n+p)}{(2n) \Gamma(2n)}\, \frac{1}{x^{2n+p}}
\label{asymptoticpbder}
\ear
This function also has the following integral representation:
\bear
\bar{\psi}^{(p)}(x)= -(-2)^{p}\int_0^{\infty}ds\,{\rm e}^{-2xs}
\,s^p\, \Bigl(\frac{1}{\sinh s} - {1\over s}\Bigr)
\label{integralpbder}
\ear
As in the Miki case,  
we can use (\ref{integralpbder}) to define $\bar{\psi}^{(p)}(x)$ for non-integer
$p$.  We can then derive new identities, for any positive $p$, by 
squaring these two representations of $\bar{\psi}^{(p)}(x)$, and then comparing, 
just as was done (for $p=0$) to prove the FPZ identity.
The proof proceeds in a very similar manner.

\underline{\sl Lemma 4.2.1 :}
\bear
\left[\bar{\psi}^{(p)}(x)\right]^2 \sim \sum_{n=2}^\infty 
\frac{1}{x^{2n+2p}} 
\sum_{k=1}^{n-1} \frac{\bar{B}_{2k} 
\bar{B}_{2n-2k}}{(2k)(2n-2k)}\frac{\Gamma
(2k+p)\Gamma(2n-2k+p)}{ \Gamma(2k)\Gamma(2n-2k)}\, 
\label{asymptoticptdersq}
\ear
\underline{\sl Proof:} follows from (\ref{asymptoticpbder}). $\square$

\underline{\sl Lemma 4.2.2 :}
\bear 
\left[\bar{\psi}^{(p)}(x)\right]^2&=& 2^{2p+1}\int_0^\infty dy\, 
y^{2p+1}\, e^
{-2 x y}\int_0^1 du \, u^p(1-u)^p\left\{\frac{1}{\sinh y}\left(\coth 
yu-\frac{1}
{yu}\right) \right.
\nonumber\\
&& \left. -\frac{1}{y(1-u)}\left[u\left( \frac{1}{\sinh 
yu}-\frac{1}{yu}
\right)- \left( \frac{1}{\sinh y}-\frac{1}{y}\right)\right] \right\}
\label{l422}
\ear
\underline{\sl Proof :} Square the integral representation (\ref
{integralpbder}), change variables from $s$ and $s^\prime$ to $y$ and 
$u$, 
and regroup terms as in the proof of {\sl Lemma 2.2.1}. 
$\square$

Now consider each of the two terms appearing on the RHS of 
(\ref{l422}).

\underline{\sl Lemma 4.2.3 :}
\bear 
&&2^{2p+1}\int_0^\infty dy\, y^{2p+1}\, e^{-2 x y}\int_0^1 du \, 
u^p(1-u)
^p\left\{\frac{1}{\sinh y}\left(\coth yu-\frac{1}{yu}\right) 
\right\}\nonumber\\
&&\quad \sim \quad 2 \Gamma(p+1) \sum_{n=1}^\infty \frac{1}{x^{2n+2p}} 
\sum_{k=1}^{n} 
\frac{B_{2k} 
\bar{B}_{2n-2k}}{(2k)!(2n-2k)!}\frac{\Gamma(2k+p)\Gamma(2n+2p)}
{\Gamma(2p+2k+1)}
\label{l423}
\ear
\underline{\sl Proof :} The proof is almost identical to the proof of 
{\sl 
Lemma 4.1.4}, but in doing the $y$ integral we use an asymptotic 
expansion of 
$1/\sinh y$ rather than $\coth y - {1\over y}$. This has the effect of replacing one 
of the 
Bernoulli number factors $B_{2n-2k}$ by $\bar{B}_{2n-2k}$, and also of changing the upper 
limit of 
the $k$ summation from $(n-1)$ to $n$. $\square$

\underline{\sl Lemma 4.2.4 :}
\bear 
&&- 2^{2p+1}\int_0^\infty dy\, y^{2p+1}\, e^{-2 x y}\int_0^1 du \, 
u^p(1-u)
^p\left\{\frac{1}{y(1-u)}\left[u\left( \frac{1}{\sinh 
yu}-\frac{1}{yu}\right)- 
\left( \frac{1}{\sinh y}-\frac{1}{y}\right)\right] \right\}\nonumber\\
&&\quad \sim \quad 2 \sum_{n=2}^\infty \frac{1}{x^{2n+2p}} 
\frac{\bar{B}_{2n} 
\Gamma(2n+2p)}{(2n)!} \sum_{k=1}^{2n} \beta(p+k,p+1)
\label{l424}
\ear
\underline{\sl Proof :} The proof is almost identical to the proof of 
{\sl 
Lemma 4.1.5}, except we use the asymptotic expansion of $1/\sinh y$ 
rather than 
$\coth y$. $\square$

We are now ready to state the generalization of the FPZ identity:

\underline{\sl Theorem 4.2 :} For any $p\geq 0$, and for 
integer 
$n\geq 2$:
\bear
\sum_{k=1}^{n-1} \frac{\bar{B}_{2k} 
\bar{B}_{2n-2k}}{(2k)(2n-2k)}\frac{\Gamma
(2k+p)\Gamma(2n-2k+p)}{ \Gamma(2k)\Gamma(2n-2k)} &=& 2\Gamma(p+1)\sum_{k=1}^{n} 
\frac{B_
{2k} 
\bar{B}_{2n-2k}}{(2k)!(2n-2k)!}\frac{\Gamma(2k+p)\Gamma(2n+2p)}{\Gamma
(2p+2k+1)}\nonumber\\
&& +2  \frac{\bar{B}_{2n} \Gamma(2n+2p)}{(2n)!} \sum_{k=1}^{2n-1} 
\beta
(p+k,p+1)
\label{t42}
\ear

\underline{\sl Proof :} Follows by comparing the result of {\sl 
Lemma 4.2.1} 
with those of {\sl Lemmas} 4.2.2 - 4.2.4. $\square$

\underline{\sl Comment 4.2.1 :} When $p=0$ we recover from {\sl 
Theorem 4.2} 
the FPZ identity (\ref{fpzidentity}). 

\underline{\sl Comment 4.2.2 :} When $p=1$ we obtain from {\sl Theorem 
4.2} a 
convolution identity which just involves the Bernoulli numbers 
themselves on 
the left-hand side: for $n\geq 1$,
\bear
\sum_{k=1}^{n} \bar{B}_{2k} \bar{B}_{2n-2k}= \frac{1}{n+1} 
\sum_{k=1}^{n} B_
{2k} \bar{B}_{2n-2k}{2n+2\choose 2k+2} + 2n \bar{B}_{2n}
\label{c422}
\ear

\subsection{Generalization of Theorem 3.1}
\label{mixedgen}

To generalize {\sl Theorem 3.1}, we differentiate $p$ 
times the relation (\ref{crossrelation}) connecting the two generating 
functions $\tilde{\psi}(x)$ and $\bar{\psi}(x)$. This leads to

\bear
 \tilde{\psi}^{(p)}(x)+ \bar{\psi}^{(p)}(x)=2^{p+1} 
\tilde{\psi}^{(p)}(2x)
\label{relmixed}
\ear
This relation also holds true for non-integer positive $p$, as can be easily
seen using the integral representations (\ref{integralptder}), (\ref{integralpbder})
for $\tilde{\psi}^{(p)}$ and $\bar{\psi}^{(p)}$ and the trigonometric identity
\bear
{\rm coth}(s) + {1\over {\rm sinh(s)}}
&=& {\rm coth}\bigl({s\over 2}\bigr)
\label{idtrigmixed}
\ear 
Squaring the relation (\ref{relmixed}) we obtain

\underline{\sl Lemma 4.3.1 :}
\bear
2 \tilde{\psi}^{(p)}(x)\, \bar{\psi}^{(p)}(x)=2^{2p+2} \left[ 
\tilde{\psi}^
{(p)}(2x)\right]^2 - \left[ \tilde{\psi}^{(p)}(x)\right]^2 - \left[ 
\bar{\psi}^
{(p)}(x)\right]^2
\label{l431}
\ear
($p\geq 0$). This brings us to

\underline{\sl Theorem 4.3 :} For any  $p\geq 0$, and for 
integer 
$n\geq 2$:
\bear
&&\sum_{k=1}^{n-1} \frac{B_{2k} 
\bar{B}_{2n-2k}}{(2k)(2n-2k)}\frac{\Gamma(2k+p)
\Gamma(2n-2k+p)}{ \Gamma(2k)\Gamma(2n-2k)} = 
\hspace{180pt}
\nonumber\\
&&\hskip 2.1cm 2\Gamma(p+1) \sum_{k=1}^{n}\frac{B_{2k} 
B_{2n-2k}}{(2k)!(2n-2k)!} \Bigl
({1-2^{2k-1}\over 2^{2n-1}}\Bigr) 
\frac{\Gamma(2k+p)\Gamma(2n+2p)}{\Gamma
(2p+2k+1)}
\nonumber\\&&\hskip 2.1cm
+\frac{B_{2n} \Gamma(2n+2p)}{(2n)! 2^{2n-1}} 
\sum_{k=1}^{2n-1} \beta
(p+k,p+1)\nonumber\\
\label{t43}
\ear

\underline{\sl Proof :} The proof follows by taking the product of the 
expansions of the two functions on the left-hand side of {\sl Lemma 
4.3.1}, and 
comparing  with the expansions of the squares of the integral 
representations 
of the three terms appearing on the right-hand side of {\sl Lemma 4.3.1}, using the results of 
{\sl 
Theorem 4.1} and {\sl Theorem 4.2}. $\square$

\underline{\sl Comment 4.3.1 :} When $p=0$ we recover the identity in {\it Theorem 3.1}.

\underline{\sl Comment 4.3.2 :} When $p=1$ we obtain from {\sl Theorem 4.3} a 
convolution identity which just involves the Bernoulli numbers themselves on 
the left-hand side: for $n\geq 1$,
\bear
\sum_{k=1}^{n-1} B_{2k} \bar{B}_{2n-2k}= \frac{1}{n+1} \sum_{k=1}^{n} 
B_{2k} B_
{2n-2k}  \Bigl({1-2^{2k-1}\over 2^{2n-1}}\Bigr) {2n+2\choose 2k+2} 
+(2n-1) \frac
{B_{2n}}{2^{2n}}
\label{c431}
\ear

\underline{\sl Comment 4.3.3 :} Note that in all the above the positiveness condition on
$p$ was used only to avoid singularities. Theorems 4.1, 4.2 and 4.3 actually hold true also
for negative $p$ as long as none of the $\Gamma$ - factors on either side becomes singular.

\section{Higher order convolution identities}
\label{multi}
\renewcommand{\theequation}{5.\arabic{equation}}
\setcounter{equation}{0}
In the recent \cite{gessel}, I.M. Gessel shows the existence of an infinite tower
of convolution identities involving multiple products of Bernoulli numbers, of which
Miki's identity (\ref{mikiidentity}) is just the lowest order one. He also explicitly obtains the
next element of this series, a triple product identity:

\underline{{\sl Theorem 5.1 :}} (\cite{gessel}, eq.(4)): For integer $n\geq 3$,
\bear
\sum_{k+l+m=n\atop k,l,m \geq 1}
{B_{2k}B_{2l}B_{2m}\over (2k)(2l)(2m)}
&=&
\sum_{k+l+m=n\atop k,l,m \geq 1}
{B_{2k}B_{2l}B_{2m}\over (2k)(2l)(2m)}
{2n \choose 2k,2l,2m}
+3H_{2n}\sum_{k=1}^{n-1}{2n\choose 2k} {B_{2k}B_{2n-2k}\over (2k)(2n-2k)}
\nonumber\\
&& +6H_{2n,2}{B_{2n}\over 2n}
-\Bigl(n^2-{3\over 2}n+{5\over 4}\Bigr) {B_{2n-2}\over (2n-2)}
\label{gesselidentity}
\ear
Here
\bear
H_{2n,2} &\equiv & \sum_{1\leq i < j \leq 2n} {1\over ij}\quad
\biggl( =\sum_{i=1}^{2n-1}{H_l\over l+1}\biggr)
\label{defH2}
\ear
In our present approach, it is clear how to construct generalizations of the
Miki and FPZ identities involving $N$ - fold products of Bernoulli numbers:

\begin{enumerate}

\item
Take the $N$th power of $\tilde\psi(x)$ (resp. $\bar\psi(x)$). 
The expansions (\ref{asymptildepsi}) (resp. (\ref{psibarlargex}))
generate the $N-1$ -- fold convolution on the left-hand side of the identity,
\bear
\sum_{n=N}^{\infty} {1\over x^{2n}}
\sum_{\sum_{i=1}^Nk_i = n\atop k_1,k_2,\ldots ,k_N \geq 1}
\,\prod_{i=1}^N{B_{2k_i}\over 2k_i}
\label{lhsmultid}
\ear
(with $B_l$ replaced by $\bar B_l$ in the FPZ case).

\item
Use (\ref{intreppsitilde}) (resp. (\ref{intreppsibar})) 
to rewrite 
\bear
\bigr(\tilde\psi(x)\bigl)^N
&=&
\int_0^{\infty}ds_1\int_0^{\infty}}ds_2
\cdots\int_0^{\infty}ds_N\, {\rm e}^{-2x(s_1+s_2+\ldots +s_N)}
\,\prod_{i=1}^N\Bigl({\rm coth(s_i)-{1\over s_i}\Bigr)
\nonumber\\
\label{rhsmultid}
\ear
(with ${\rm coth}$ replaced by ${1\over \sinh}$ in the FPZ case).

\item
Use trigonometric identities together with the symmetry of the integrand under
permutations of the variables $\lbrace s_1,\ldots,s_N \rbrace$
to rewrite the integrand in such a way that only 
the combinations
$s_1+s_2+\ldots +s_N$, $s_2+s_3+\ldots +s_N$, \ldots, $s_{N-1}+s_N$, $s_N$
appear as arguments of trigonometric functions (it is easy to see that this is always
possible).

\item
Perform the transformation of variables from $\lbrace s_1,\ldots,s_N \rbrace$ to 
$\lbrace y,u_1,u_2,\ldots,u_{N-1}\rbrace$ where
\bear
y &=& s_1+s_2+\ldots + s_N \nonumber\\
u_M &=& {s_{M+1}+s_{M+2}+\ldots + s_N \over s_1+s_2+\ldots + s_N},\quad 
M=1,\ldots,N-1
\label{trafogeneral}
\ear
The Jacobi factor of this transformation is $y^{N-1}$.

\item
Use the Taylor expansions (\ref{taylcoth}),(\ref{taylsinh}) to do all integrals. 

\end{enumerate}
Let us carry this through explicitly for the case $N=3$. In the Miki case, after step
3 one finds

\noindent
\underline{{\sl Lemma 5.1:}}
\bear
\bigr(\tilde\psi(x)\bigl)^3
&=&
\int_0^{\infty}ds_1\int_0^{\infty}ds_2\int_0^{\infty}ds_3\, {\rm e}^{-2x(s_1+s_2+s_3)}
\nonumber\\&&\times
\biggl\lbrace
6C_{123}C_{23}C_3
+ {6\over s_1}\Bigl\lbrack C_{123}-{s_{23}\over s_{123}}C_{23}\Bigr\rbrack C_3
+{6\over s_2}C_{123}\Bigl\lbrack C_{23}-{s_3\over s_{23}}C_3\Bigr\rbrack
\nonumber\\&&
+{6\over s_2}\Bigl\lbrack {1\over s_1}(C_{123}-{s_{23}\over s_{123}}C_{23})
-{1\over s_{12}}(C_{123}-{s_3\over s_{123}}C_3)\Bigr\rbrack
-3C_3-2C_{123}-{4\over s_{123}}
\biggr\rbrace\nonumber\\
\label{intmikitriple}
\ear 
Here for compactness we have introduced the abbreviations
\bear
s_{23}= s_2+s_3, \quad s_{123} = s_1+s_2+s_3,\quad
C_{(\cdot)}=\coth(s_{(\cdot)})-{1\over s_{(\cdot)}}.
\label{abbr}
\ear
Moreover, we have already combined terms in a way which facilitates the evaluation
of the integrals (in particular, it avoids the appearance of spurious singularities).

After the transformation (\ref{trafogeneral}) the integrals can be done in a way which
is completely analogous to the $N=2$ case treated in section \ref{mikifpz}. The result
is a slight modified form of the identity (\ref{gesselidentity}),

\underline{{\sl Theorem 5.2 :}}
(Modified form of Gessel's identity): For 
integer $n\geq 3$, 
\bear
\sum_{k+l+m=n\atop k,l,m \geq 1}
{B_{2k}B_{2l}B_{2m}\over (2k)(2l)(2m)}
&=&
{3\over 2n}\sum_{k+l+m=n\atop k,l,m \geq 1}
{B_{2k}B_{2l}B_{2m}\over (2k)(2l)}
{2n \choose 2k,2l,2m}
+{3\over n}H_{2n}\sum_{k=1}^{n-1}{2n\choose 2k} {B_{2k}B_{2n-2k}\over (2k)}
\nonumber\\
&& +6H_{2n,2}{B_{2n}\over 2n}
-\Bigl(n^2-{3\over 2}n+{5\over 4}\Bigr) {B_{2n-2}\over (2n-2)}
\label{gesselidentitymod}
\ear

\noindent
\underline{{\sl Comment 5.1 :}}
The somewhat different form of the right hand side compared to (\ref{gesselidentity}) is due
to the same type of ambiguity (regarding integrations-by-parts) which was mentioned already at the end of section 2.1. 

The FPZ case is again similar, though slightly more complicated. After step three one finds that the integrand can be written in the following way,

\noindent
\underline{{\sl Lemma 5.2:}}
\bear
\bigr(\bar\psi(x)\bigl)^3
&=&
\int_0^{\infty}ds_1\int_0^{\infty}ds_2
\int_0^{\infty}ds_3\, {\rm e}^{-2x(s_1+s_2+s_3)}
\nonumber\\&&\times
\biggl\lbrace
6S_{123}C_{23}C_3
+ {6\over s_1}\Bigl\lbrack S_{123}-{s_{23}\over s_{123}}S_{23}\Bigr\rbrack C_3
+{6\over s_2}S_{123}\Bigl\lbrack C_{23}-{s_3\over s_{23}}C_3\Bigr\rbrack
\nonumber\\&&
+{6\over s_2}\Bigl\lbrack {1\over s_1}(S_{123}-{s_{23}\over s_{123}}S_{23})
-{1\over s_{12}}(S_{123}-{s_3\over s_{123}}S_3)\Bigr\rbrack
-2S_{123}-{4\over s_{123}}
\nonumber\\&&
+{6\over s_{123}}[C_{23}-S_{23}]C_3 +{6\over s_2s_{123}}\lbrack C_{23}-S_{23}-C_3+S_3\rbrack
\biggr\rbrace\nonumber\\
\label{intfpztriple}
\ear 
Here we have used one more abbreviation,
\bear
S_{(\cdot )} &=& {1\over \sinh(s_{(\cdot)})} - {1\over s_{(\cdot)}}
\label{abbr2}
\ear

\underline{{\sl Theorem 5.3 :}}
(Cubic generalization of the FPZ identity): For 
integer 
$n\geq 3$, 
\bear
\sum_{k+l+m=n\atop k,l,m \geq 1}
{\bar B_{2k}\bar B_{2l}\bar B_{2m}\over (2k)(2l)(2m)}
&=&
{3\over 2n}\sum_{k+l+m=n\atop k,l,m \geq 1}
{B_{2k}B_{2l}\bar B_{2m}\over (2k)(2l)}
{2n \choose 2k,2l,2m}
+{3\over n}H_{2n}\sum_{k=1}^{n-1}{2n\choose 2k} {B_{2k}\bar B_{2n-2k}\over (2k)}
\nonumber\\
&& +{3\over 2n^2}\sum_{k=1}^{n-1}{2n\choose 2k}{B_{2k}\over 2k}(B_{2n-2k}-\bar B_{2n-2k})
+{3\over 2n^2}H_{2n-1}(B_{2n}-\bar B_{2n})
\nonumber\\&&
+6H_{2n,2}{\bar B_{2n}\over 2n}
-{2n-1\over 4}\bar B_{2n-2}
\label{fpzcubic}
\ear
Note again the similarity with the Miki case, (\ref{gesselidentitymod}).

\noindent
\underline{{\sl Comment 5.2 :}}
It would be straightforward to extend both (\ref{gesselidentitymod}) and
(\ref{fpzcubic}) to continuous families of identities along the lines of
section 4.

\section{Conclusions}
\label{conclusions}
\renewcommand{\theequation}{6.\arabic{equation}}
\setcounter{equation}{0}

The method presented here allows one to derive, with relatively little effort, 
convolution identities for Bernoulli numbers of the quadratic type
as well as higher order ones. Clearly we have not been able here to
explore all its ramifications. For example, it should be possible to derive
``mixed'' identities such as (\ref{t31}) also at the cubic or higher level.
Another possible direction is to use other generating functions to 
generate 
related identities involving the Euler numbers. The simplest such case comes 
from considering the generating function
\bear
g(x)&=&\int_0^\infty ds\, e^{-2xs}\, {\rm sech}\, s \nonumber
\\
&\sim&\sum_{n=0}^\infty \frac{E_{2n}}{(2x)^{2n+1}}
\label{eulermiki}
\ear
Then it follows that 
\bear
\left[g(x)\right]^2&=&2 \int_0^\infty dy\, e^{-2xy}\, \frac{\ln \cosh 
y}{\sinh 
y} \quad ,
\label{eulermikisq}
\ear
from which one finds
\bear
\sum_{k=1}^n E_{2k-2} E_{2n-2k} =\frac{2}{n} \sum_{k=1}^n \frac{B_{2k} 
B_{2n-
2k} }{k} \left(2^{2k} - 1\right)2^{2k - 1}\left(1 - 2^{2n - 2k - 
1}\right) {2n 
\choose 2k}
\label{eulernew}
\ear
($n\geq 1$).
Clearly one can generate many other such identities relating 
convolutions of 
Euler numbers to convolutions of Bernoulli numbers.

We conclude by emphasizing again that the types of generating functions 
and 
identities discussed here show up naturally in perturbative quantum 
field 
theory computations at the second-order (or ``two loop") level \cite{sd1,sd2,loops}.  We 
expect 
related multiple convolution identities of higher order to play a similar role for higher-loop 
contributions 
to the effective Lagrangian in quantum electrodynamics beyond the 
two-loop 
level. 
The higher order FPZ type identities might correspond to new relations
between Hodge integrals and thus be of relevance for topological
quantum field theory and string theory.

{\bf Note added (January 2013):} This article was submitted to the 
arXiv in June 2004, as \url{http://arxiv.org/abs/math/0406610} but was delayed in publication review. Since then some related work citing our preprint has appeared:
In \cite{crabb} a method based on generating functions, similar to the one
introduced here, was used to derive a convolution identity for Bernoulli polynomials that generalizes both Miki's identity (\ref{mikiidentity}) and the FPZ identity (\ref{fpzidentity}). It was also outlined how to use the same method to generalize Theorem 4.1 to the (univariate) polynomial level. In \cite{pansun2005} this convolution identity was further generalized to
a bivariate one. In the same article, it was noted that our (\ref{c412}) had already been stated as a conjecture by Matiyasevich in 1997 \cite{matisayevich}, and a bivariate polynomial generalization of this identity was obtained, too. 
In \cite{fupanzhang} Theorem 4.1 was further generalized to the bivariate polynomial level, and the bivariate quadratic
convolution identity of \cite{pansun2005} was generalized to a family of multivariate multiple convolution identities for Bernoulli polynomials at any order (at the cubic level those presumably also relate to the identities of section 5, although no such claim was made in \cite{fupanzhang}). In \cite{gadpad} a modification of our generating function $\tilde \psi$
was pointed out that may lead to another type of generalization of Miki's identity. Finally, in \cite{gorzhi} the
known vanishing of a certain type of one-loop amplitudes in N=4 Super-Yang-Mills theory was used to derive yet
another type of quadratic convolution identities involving Bernoulli numbers, not obviously related to any of the
above. These results  further strengthen the case for an ubiquitous role of Bernoulli number identities in
perturbative quantum field theory.

{\bf Acknowledgements:} We are very grateful to Richard Stanley for 
correspondence, and to Albrecht Klemm for discussions.
C.S. thanks the Institut des Hautes \'Etudes 
Scientifiques,
Bures-sur-Yvette, and the Albert-Einstein-Institut, Potsdam, for hospitality. 
G.D. thanks the US Department of 
Energy for 
support under grant DE-FG02-92ER40716, and thanks the Rockefeller 
Foundation 
for a Bellagio Residency Award. We also acknowledge the support of NSF 
and 
CONACyT for a US-Mexico collaborative research grant NSF-INT-0122615.

\end{document}